\documentclass[11pt,reqno]{amsart}
\usepackage[T1]{fontenc}
\usepackage[utf8]{inputenc}
\usepackage{amsmath,amssymb,amsfonts,amsthm,mathrsfs,mathtools}
\usepackage[margin=1.5in]{geometry}
\usepackage{graphicx}
\usepackage{url}
\allowdisplaybreaks

\theoremstyle{plain}
\newtheorem{theorem}{Theorem}[section]
\newtheorem{lemma}[theorem]{Lemma}

\theoremstyle{definition}

\theoremstyle{remark}
\newtheorem{remark}[theorem]{Remark}
\newtheorem*{remark*}{Remark}
\numberwithin{equation}{section}
\newcommand{\RNum}[1]{\uppercase\expandafter{\romannumeral #1\relax}}

\title[Characterizations of the plane and the catenoid]{Characterizations of the plane and the catenoid as capillary minimal surfaces}
\author[E. Yeon]{Eungbeom Yeon}
\address{Department of Mathematical Sciences, Seoul National University, Seoul 08826, Republic of Korea}
\curraddr{Department of Mathematics and Institute of Mathematical Science, Pusan National University, Busan 46241, Republic of Korea}
\email{y2bum@snu.ac.kr}
\date{}

\begin{document}

\begin{abstract}
In this paper, we prove that an immersed annular capillary minimal surface outside the unit ball in $\mathbb{R}^3$ with one complete embedded end and finite total curvature must be either part of the plane or part of the catenoid. The proofs are based on the use of the Weierstrass data through the Hopf differential and associated holomorphic functions. For a planar end, we control the second fundamental form, while for a catenoidal end, we analyze the Gauss map in more detail using the boundary condition. We also prove nonexistence results for annular free boundary minimal surfaces with one end asymptotic to the end of Enneper's surface, under appropriate conditions. Finally, we establish an exterior analogue of Choe's result on capillary minimal surfaces with vertices.
\end{abstract}

\maketitle

\section{Introduction}\label{sec:intro}

The plane and the catenoid are two important examples in minimal surface theory. The plane is the only totally umbilic minimal surface in $\mathbb{R}^3$ and the catenoid is the only nonplanar rotationally symmetric minimal surface in $\mathbb{R}^3$. Besides these properties, various characterizations of the plane and the catenoid have been studied for many years. For example, the plane is the only entire minimal graph in $\mathbb{R}^3$ \cite{B}. A complete minimal surface of finite total curvature with two embedded ends must be the catenoid \cite{S}. In fact, the catenoid is the only embedded complete minimal annulus with finite total curvature \cite{JM}. It is also known that a complete minimal annulus in $\mathbb{R}^3$ whose intersection with every $z$-plane is a Jordan curve must be the catenoid \cite{Col}. In \cite{BB}, it was proved that among properly embedded minimal annuli spanning the two boundary planes of a slab, the maximally stable catenoidal waist symmetric with respect to the middle plane has the minimum area. After a suitable scaling and rigid motion, this waist becomes the critical catenoid in the unit ball.

A capillary minimal surface $\Sigma$ in a domain $U$ is a minimal surface which meets $\partial U$ at a constant contact angle in $(0,\pi)$ along $\partial\Sigma$. When the contact angle is $90^\circ$, the surface is called a free boundary minimal surface. Compared with the free boundary case, results on capillary angles are more limited. A classical result of Nitsche \cite{N} states that a disk type capillary minimal surface inside the unit ball must be a planar disk. The proof uses the Hopf differential together with the fact that the boundary is a line of curvature. Ros and Souam \cite{RS} generalized this result to balls in the space forms $\mathbb{H}^3$ and $\mathbb{S}^3$. Fern\'andez, Hauswirth and Mira \cite{FHM} constructed non-rotational free boundary minimal annuli immersed in the unit ball, giving a negative answer to Nitsche's question on the uniqueness of the immersed annuli. For capillary annuli with planar boundaries, Pyo \cite{P} showed that an immersed minimal annulus with two planar Jordan boundary curves lying in distinct planes and meeting those planes at constant contact angles must be part of the catenoid.

In this paper, we first characterize the plane and the catenoid as capillary minimal surfaces with embedded ends. Instead of looking at the capillary minimal surface inside the unit ball $\mathrm{B}^3$, we consider the exterior problem that generalizes Nitsche's theorem \cite{N}. More specifically, we study a capillary minimal surface in $\mathbb{R}^3 \setminus \mathrm{B}^3$ with one embedded end and with boundary lying on the unit sphere $\mathrm{S}^2$. To do this, we use two classical theorems of Osserman and Schoen. Osserman \cite{O} proved that a complete minimal surface with finite total curvature is conformally equivalent to a compact Riemann surface with a finite number of punctures. An end corresponding to a puncture is complete if every curve tending to the puncture has infinite length in the induced metric. Schoen \cite{S} showed that an embedded end of a complete minimal surface with finite total curvature must be a planar or a catenoidal one. 

Let
\begin{align*}
\mathrm D &= \{z\in\mathbb C \mid |z|<1\},\\
\mathrm D^\prime &= \{z\in\mathbb C \mid 0<|z|<1\},\\
\overline{\mathrm D}^\prime &= \{z\in\mathbb C \mid 0<|z|\le1\}.
\end{align*}
As in Nitsche's theorem, we work with conformally parametrized minimal surfaces. In our exterior setting, we consider a conformal immersion $X:\mathrm D^\prime\to\mathbb R^3$ whose image $\Sigma$ has finite total curvature and one embedded end. Here the puncture corresponds to the end of $\Sigma$, and $\partial\mathrm D$ to its boundary $\partial\Sigma$. The end is called embedded if the restriction of $X$ to $\{z\in\mathbb C\mid 0<|z|<\varepsilon\}$ is an embedding for some $0<\varepsilon<1$. We have the following theorem, which gives an exterior version of Nitsche's theorem \cite{N} under the assumptions of finite total curvature and one embedded end.
\begin{theorem}\label{thm1}
Let $X\in C^2(\mathrm{D}^\prime,\mathbb R^3)\cap C^1(\overline{\mathrm{D}}^\prime,\mathbb R^3)$ be a capillary minimal immersion in $\mathbb{R}^3\setminus \mathrm{B}^3$ with finite total curvature and one complete embedded end. Assume that $X$ is a regular conformal immersion up to the boundary and is $C^{2,\gamma}$ in a neighborhood of $\partial\mathrm D$ for some $0<\gamma<1$. Then its image must be either part of the plane or part of the catenoid.
\end{theorem}
The proof of Theorem~\ref{thm1} makes effective use of the Weierstrass data. In the case of a planar end, control of the second fundamental form obtained from the capillary condition and the asymptotic behavior at the end directly leads to the result. In the case of a catenoidal end, we also need a more detailed analysis of the Gauss map using the boundary condition for the induced metric. Although boundary conditions can be difficult to express in terms of the Weierstrass data, the capillary condition on the sphere gives useful relations for the second fundamental form and the metric. In the catenoidal case, the boundary condition and the behavior at the end give inequalities involving the Gauss map. The essential part of the proof is to show that equality holds in these inequalities. This argument is inspired by Osserman and Schiffer's analysis of doubly-connected minimal surfaces \cite{OS}. In this regard, the Weierstrass data play an important role.

Enneper's surface is a simply connected, complete immersed minimal surface in $\mathbb{R}^3$ with one immersed end. It is well known that the only complete orientable minimal surfaces with total curvature $-4 \pi$ are the catenoid and Enneper's surface \cite{O}. So it can be said that Enneper's surface is the simplest complete minimal surface in $\mathbb{R}^3$ whose end is not embedded. The end of Enneper's surface has multiplicity 3. We will say that a surface has an end asymptotic to the end of Enneper's surface if there are a local conformal coordinate $z$ with the puncture at $z=0$ and a conformal parametrization $X_E$, on the same punctured neighborhood, of an end of Enneper's surface up to a rigid motion and scaling, such that
\begin{align*}
X(z)-X_E(z)=O(1)\qquad\text{as }z\to0.
\end{align*}
Since $X-X_E$ is bounded and harmonic, it extends harmonically across the puncture. Thus after a translation of $X_E$, the distance between the two ends tends to zero. 

It is natural to ask whether a result similar to Theorem~\ref{thm1} holds when an embedded end is replaced by an end asymptotic to the end of Enneper's surface. It turns out that more complicated cases arise for Enneper ends, since the Hopf differential gives less control over the second fundamental form than in the embedded end case. We obtain the following nonexistence result for free boundary minimal surfaces under some additional assumptions.
\begin{theorem}\label{thm2}
Let $X\in C^2(\mathrm D^\prime,\mathbb R^3)\cap C^1(\overline{\mathrm D}^\prime,\mathbb R^3)$ be a conformal free boundary minimal immersion in $\mathbb R^3\setminus\mathrm B^3$ which is regular up to the boundary and is $C^{2,\gamma}$ near $\partial\mathrm D$ for some $0<\gamma<1$. Assume that the surface has finite total curvature and one complete end asymptotic to the end of Enneper's surface. If its boundary is embedded and weakly convex in $\mathrm S^2$, then such an immersion does not exist.
\end{theorem}
The main idea of the proof is to combine Hopf's methods \cite{choe,H1} with the geometry of the boundary. For an Enneper end, more complicated behavior of the curvature lines at boundary umbilic points is not excluded. We therefore assume that the boundary is embedded and weakly convex in $\mathrm S^2$. Here weak convexity means that its geodesic curvature is nonnegative with respect to the inward unit conormal of one of the regions it bounds. Together with embeddedness, this implies that the boundary is spherically convex and lies in a closed hemisphere \cite[Lemma~2.2]{Ghomi}. We expect that counterexamples exist when the weak convexity assumption is removed in Theorem~\ref{thm2}; see Remark~\ref{rmk:enneper}.

In \cite{choe}, characterizations of constant mean curvature capillary surfaces were studied by generalizing Hopf's methods. This generalized Nitsche's theorem to domains bounded by several spheres or planes. Motivated by this work, we generalize Theorem~\ref{thm1} to the exterior of convex domains bounded by spheres and planes. Each domain is understood as the common interior of the corresponding balls and half-spaces. Here each regular boundary arc lies in a single face, and the capillary contact angles may differ between faces. As in \cite{choe}, we call the singular points of the piecewise smooth boundary of the surface its vertices. If the boundary has at least one singular vertex, we obtain the following theorem.

\begin{theorem}\label{thm3}
Let $U$ be a convex domain in $\mathbb R^3$ bounded by finitely many spheres and planes. Let $X:\overline{\mathrm D}^\prime\to\mathbb R^3$ be a conformal parametrization of a capillary minimal surface $\Sigma$ whose interior is contained in $\mathbb R^3\setminus\overline U$ and whose boundary lies on $\partial U$. Suppose that $\Sigma$ has finite total curvature and one complete embedded end. Assume that $\Sigma$ is $C^{2,\gamma}$ and regularly immersed up to and including $\partial\Sigma$, for some $0<\gamma<1$. Its boundary is piecewise $C^{2,\gamma}$ with at least one vertex, and the two boundary tangent lines at each vertex are distinct. Then $\Sigma$ cannot have a catenoidal end. Furthermore, if its end is planar, $\Sigma$ must be part of the plane.
\end{theorem}

This paper is organized as follows. In Section~\ref{sec2}, we review some of the important facts which will be used throughout the rest of the paper. In Section~\ref{sec3}, we prove Theorem~\ref{thm1}. In Section~\ref{sec4}, we prove Theorem~\ref{thm2}. Finally, in Section~\ref{sec5}, we prove Theorem~\ref{thm3}.

\section*{Acknowledgement}
Problems in this paper were proposed to me by my advisor, Jaigyoung Choe, in the fall of 2017 for my doctoral dissertation. I would like to express my gratitude to him for his fruitful advice and guidance throughout my doctoral studies. A few modifications to the proof of Theorem~\ref{thm1} were made in response to comments received during the peer review process. I thank the anonymous reviewer for the helpful comments.

\section{Preliminaries}\label{sec2}

In this section we present the Weierstrass representation formula for minimal surfaces in $\mathbb{R}^3$. We know that the Gauss map of a complete minimal surface of finite total curvature $\Sigma$ can be considered as a meromorphic function $g(z) : \Sigma \rightarrow \mathbb{C}\cup\{\infty\}$ where $z=u+iv$ is a local conformal coordinate on the surface. Taking the third coordinate function $x_3$ defined on the surface, we have the holomorphic differential $dh = dx_3 + i dx_3^*$. Here, we have to be aware that the harmonic conjugate $x_3^*$ can be a multi-valued function since the surface need not be simply connected. The pair $(g,dh)$ is called the Weierstrass data and the minimal surface $X : \Sigma \rightarrow \mathbb{R}^3$ can be represented by the real part of the holomorphic curve as follows.
\begin{equation}\label{eq1}
X(p) = \text{Re} \int_{p_0}^p \left(\frac{1}{2}\left(\frac{1}{g}-g\right) ,\frac{i}{2}\left(\frac{1}{g}+g\right) , 1\right) dh.
\end{equation}

Conversely, if the Weierstrass data $(g,dh)$ are given, the representation~\eqref{eq1} gives rise to a conformal minimal immersion if the zeroes of $dh$ coincide with the zeroes and poles of $g$ with the same order, and for any closed curve $\gamma$ in $\Sigma$,
\begin{align*}
\overline{\int_\gamma g dh} = \int_\gamma \frac{dh}{g} , \quad \text{Re}\int_\gamma dh = 0
\end{align*}
holds.

With the Weierstrass data, various geometric invariants can also be represented by the data. In particular, the first and second fundamental forms of the surface are given by
\begin{align}\label{eq2}
ds^2 &= \left(\frac{1}{2}(|g|+|g|^{-1})|dh| \right)^2 =  \Lambda^2 |dz|^2,\\
b(v,v) &= -\text{Re} \left( \frac{dg}{g}(v)\cdot dh(v) \right),\label{eq3}
\end{align}
where $v$ is a tangent vector to the surface $\Sigma$ and $N=(X_u\times X_v)/\Lambda^2$ is the chosen unit normal. Also, the Gaussian curvature of the surface is given by
\begin{align}\label{eq4}
K = -\left(\frac{4\left|\frac{dg}{g}\right|}{(|g|+|g|^{-1})^2|dh|}\right)^2.
\end{align}
For a complete embedded end of finite total curvature, after rotating the limiting normal to the vertical direction, we have
\begin{align*}
dh=\left(\frac{a}{z}+O(1)\right)dz,\qquad x_3(z)=a\log|z|+b+O(|z|),
\end{align*}
where $a,b\in\mathbb R$ and $z$ is a local conformal coordinate with the puncture at $z=0$. By Schoen's result \cite{S}, the end is catenoidal if $a\ne0$ and planar if $a=0$.

Enneper's surface has the following Weierstrass data defined on $\mathbb C$.
\begin{align*}
g(z)=z,\qquad dh=z\,dz.
\end{align*}
The Chen-Gackstatter surface \cite{CG} is obtained by adding a handle to Enneper's surface. Up to normalization, its Weierstrass data on the elliptic curve $y^2=x(x^2-1)$ are $g=cx/y$ and $dh=dx$, where $c>0$ is chosen to satisfy the period condition. At the end, take the local coordinate $z=x^{-1/2}$. Then
\begin{align*}
g(z)=cz(1-z^4)^{-1/2}=cz+O(z^5),\qquad dh=-2z^{-3}\,dz.
\end{align*}
Recall that, in a local conformal coordinate centered at the puncture and after a translation of $X_E$, the asymptotic condition in the introduction can be written as
\begin{equation}\label{eq5}
X(z)-X_E(z)=o(1)\qquad\text{as }z\to0.
\end{equation}
The data $g_E=cz$ and $dh_E=-2z^{-3}\,dz$ represent an end of Enneper's surface up to a rotation and scaling. Since $g-g_E=O(z^5)$ and $g^{-1}-g_E^{-1}=O(z^3)$, the differences of the three differentials in~\eqref{eq1} extend holomorphically across $z=0$. Hence, after a translation of $X_E$, we have $X(z)-X_E(z)=O(|z|)$, so~\eqref{eq5} holds.

As in \cite{LR}, we can also consider the following Weierstrass data defined on $\mathbb{C}\setminus\{0\}$, where $R>1$.
\begin{equation}\label{eq6}
g(z)=z\frac{z^2-R^2}{R^2z^2-1},\qquad dh=\left(1-\frac{z^2+z^{-2}}{R^2+R^{-2}}\right)\frac{dz}{z}.
\end{equation}
The Weierstrass data in~\eqref{eq6} represent a complete minimal annulus with two ends of multiplicity 3 joined by a catenoidal neck. Indeed, the height differential in~\eqref{eq6} shows that, despite its multiplicity 3, the end at $z=0$ does not satisfy the asymptotic condition~\eqref{eq5}. Its nonzero residue at the origin gives a nonzero flux vector \cite[Section~2.3.10]{K}.

\begin{figure}[htbp]
\centering
\includegraphics[width=0.4\linewidth]{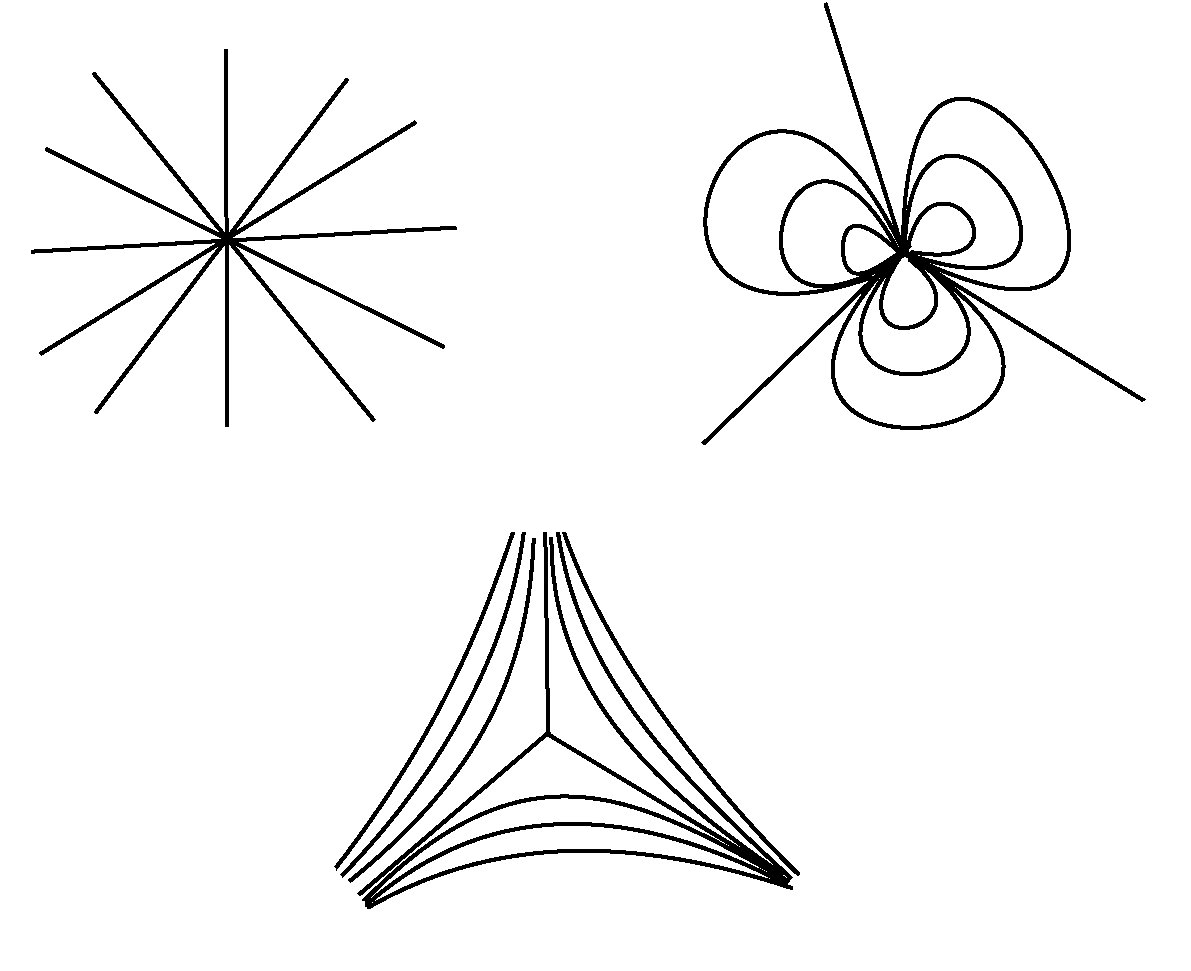}
\caption{Examples of fields of line elements with isolated singularities of indices $1$, $\frac{5}{2}$, and $-\frac{1}{2}$}
\label{fig:line-field-indices}
\end{figure}

Now, we review some facts in \cite{H1} for completeness. A field of line elements, or a line field, assigns an unoriented tangent line to each regular point of a two-dimensional surface. We assume that all fields of line elements in this paper are smooth. In case a field of line elements cannot be extended to a single point $p$, then the point $p$ is called a singularity of the field of line elements. The index of an isolated singularity $p$ of the field of line elements is defined as  $j$  where $2\pi j$ is the total change of the direction angle along a small positively oriented loop around a singular point $p$. We can see that $j$ is of the form $\frac {n}{2}$ where $n$ is an integer. The index number does not depend on what curve we choose or what metric we use on the surface. For an arbitrary integer $n$, we can consider a field of line elements with an isolated singularity of index $\frac{n}{2}$. Figure~\ref{fig:line-field-indices} shows three examples.
For a line field on a closed, orientable surface $S$ of genus $g$ with finitely many isolated singularities $p_i$, $i\in I$, the Poincar\'e--Hopf theorem \cite{H1} gives
\begin{align}
\sum_{i\in I} j_{p_i}=2-2g,\label{eq7}
\end{align}
where $j_{p_i}$ denotes the index at $p_i$.

\section{Proof of Theorem~\ref{thm1}}\label{sec3}
In this section, we prove Theorem~\ref{thm1}. As mentioned in Section~\ref{sec:intro}, we use the conformal parameter $z=u+iv$ on $\mathrm D^\prime$ and write $X:\mathrm D^\prime\to\Sigma$. Since the surface $\Sigma$ has finite total curvature, we know the fact that the Gauss map of the surface can be meromorphically extended to the puncture \cite{O}. Indeed, we see from \cite{S} that the surface must have a catenoidal end or a planar end. After rotation, we can assume that the Gauss map points at the south pole so the meromorphic function $g(z) : \mathrm D^\prime \rightarrow \mathbb {C}\cup\{\infty\}$ has the value $0$ at $z=0$. Let us denote by
\begin{align}\label{eq8}
\begin{pmatrix} \mathcal L & \mathcal M\\ \mathcal M & \mathcal N \end{pmatrix}
\end{align}
the second fundamental form of the surface as usual. Here,
\begin{align*}
\mathcal L=-\langle N_u,X_u\rangle,\quad \mathcal M=-\langle N_v,X_u\rangle,\quad \mathcal N=-\langle N_v,X_v\rangle,
\end{align*}
where $N : \Sigma \rightarrow \mathrm{S}^2$ is the Gauss map of the surface. We write $g(z)=\Pi(N(X(z)))$, where $\Pi$ denotes stereographic projection from the north pole. In terms of the above notation, the holomorphic function $\Phi$ can be written as
\begin{equation}\label{eq9}
\Phi(z)=\frac{\mathcal L-\mathcal N}{2}-i\mathcal M.
\end{equation}
As $\Phi(z)$ is known to be holomorphic on minimal surfaces \cite{H1}, the Hopf differential $\Phi (z) dz^2$ defined on minimal surfaces is very useful for dealing with curvature lines on minimal surfaces. Following Nitsche's method as in \cite{N}, we denote the real and imaginary parts of the function $f(z)=z^2\Phi(z)$ by $\alpha$ and $\beta$, respectively. Estimating the order of the Gaussian curvature in terms of the complex coordinate on the surface near the end of the surface, we have the following lemma.

\begin{lemma}\label{lem1}
Let a minimal surface $\Sigma$ be as in Theorem~\ref{thm1}, and let $\alpha+i\beta=z^2\Phi(z)$, where $\Phi$ is defined in~\eqref{eq9}. Then the harmonic functions $\alpha,\beta$ satisfy $\alpha=\beta=0$ on the whole surface when the surface has the planar end. In case the surface has a catenoidal end, $\alpha=A$ and $\beta=0$ for some nonzero constant $A$. In particular, in the catenoidal case $\Sigma$ is foliated by lines of curvature.
\end{lemma}

\begin{proof}
The Terquem-Joachimsthal theorem \cite{Spi} shows that the boundary meeting $\partial \mathrm{B}^3$ at a constant contact angle implies that $\partial \Sigma$ is a line of curvature. When we put $z=\rho e^{i\theta}$, we get by direct calculation that
\begin{equation}\label{eq10}
\frac{1}{\rho} N_{\theta} = \frac{\beta}{\rho^2 \Lambda^2} X_{\rho} + \frac{\alpha}{\rho^3\Lambda^2}X_\theta,
\end{equation}
where $N$ is the Gauss map of the surface and $\Lambda$ is defined as in \eqref{eq2}. Since the boundary of the surface becomes a line of curvature, we have   $\beta =0 $ along the boundary of the surface. We now describe the surface in terms of the Weierstrass data. If the surface is totally umbilic, the conclusion follows immediately. Otherwise, near the puncture the Weierstrass data of the surface are given by
\begin{align} \label{eq11}
g(z) &= c_nz^n+c_{n+1}z^{n+1}+\cdots,\qquad c_n\ne0,\\
\frac{dh}{g} &=\left(\frac{d_{-2}}{z^2} + \frac{d_{-1}}{z} + d_0 + d_1 z +  \cdots\right)dz,\qquad d_{-2}\ne0.\label{eq12}
\end{align}
The reason that the Gauss map can be represented as a holomorphic function near the puncture is that we assumed the Gauss map points at the south pole there. Note that the surface has a planar end for $n >1$ and a catenoidal end when $n=1$. Representation of $dh$ follows from the growth of the height function, $x_3$, which is given by
\begin{align*}
x_3(z) = \text{Re} \int dh
\end{align*}
in the Weierstrass representation formula~\eqref{eq1}.

Let $\kappa_1$ and $\kappa_2$ be the principal curvatures of the surface. Since $\begin{pmatrix} \mathcal L & \mathcal M\\ \mathcal M & \mathcal N \end{pmatrix}$ is similar to the matrix $\begin{pmatrix} \Lambda^2 \kappa_1& 0\\  0 & \Lambda^2 \kappa_2 \end{pmatrix}$, it follows that
\begin{align*}
|f(z)|^2 = \left | z^2\Phi(z) \right |^2 \le C|z|^4 \Lambda^4 \left(\kappa_1^2 +\kappa_2^2 \right) \le 2C|z|^4 \Lambda^4 |K|
\end{align*}
for some fixed constant $C$. Near the puncture, by \eqref{eq2}, \eqref{eq4}, \eqref{eq11} and \eqref{eq12}, we have
\begin{align}\label{eq13}
|f(z)|^2 \le 2C|z|^4 \Lambda^4 |K| \le C_1|z|^4|z|^{-8}|z|^{2n+2} = C_1 |z|^{2n-2}
\end{align}
near $z=0$. When the surface has a planar end, that is, $n>1$, \eqref{eq13} implies that $|f(z)|^2 \rightarrow 0 $ as $z \rightarrow 0$.  For the case $n=1$, we see that  $|f(z)|^2$ is bounded as  $z \rightarrow 0$. This means that the harmonic functions $\alpha$ and $\beta$ are bounded in a punctured disk so that $\alpha , \beta $ have removable singularities at the puncture. As a result, for arbitrary $n \ge 1$, we get $\beta \equiv 0$. Furthermore, we see that $\alpha \equiv 0 $ when the end of the surface is planar,  and $\alpha \equiv A$ in case the surface has a catenoidal end   for some nonzero constant $A$.

In the catenoidal case, $z^2\Phi=A\ne0$, so the surface has no umbilic points. Furthermore, the function $\beta$ being zero on the surface implies that \eqref{eq10} becomes
\begin{align*}
\frac{1}{\rho} N_{\theta} = \frac{\alpha}{\rho^3 \Lambda^2}X_\theta.
\end{align*}
This means that for all $p \in \mathrm D^\prime$, the vector $N_\theta \big{|}_p$ is parallel to the vector $X_\theta\big{|}_p$, which implies that every circle centered at the origin becomes a line of curvature. This completes the proof of the lemma.
\end{proof}
\begin{remark}
Fraser and Schoen \cite{FS2} used a similar holomorphic function to extend the free boundary case of Nitsche's disk theorem to higher codimension in space forms. Their function is obtained from the normal-bundle-valued Hopf differential and, in codimension one, is a constant multiple of $(z^2\Phi)^2$.
\end{remark}
We now give the proof of Theorem~\ref{thm1}.
\begin{proof}[Proof of Theorem~\ref{thm1}]
By Lemma~\ref{lem1}, we know that the harmonic functions $\alpha$ and $\beta$ both vanish throughout the whole surface when the end is planar. This implies that the surface is totally umbilic. By the condition that the surface is minimal, the surface must be a part of the plane.

Now assume that the surface has a catenoidal end. Lemma~\ref{lem1} shows that $f=z^2\Phi=A$ for some nonzero real constant $A$. In terms of the Weierstrass data, we have
\begin{equation}\label{eq14}
\Phi(z)\,dz^2=-\frac{dg}{g}\,dh=\frac{A}{z^2}\,dz^2, \qquad dh=-A\frac{g}{z^2g'}\,dz.
\end{equation}
Since the surface has no umbilic points, the Gauss map has no branch point on the punctured disk. Its simple zero at the puncture shows that it has no branch point on the closed disk as well. Here a pole of $g$ is considered in the coordinate $1/g$.

Let $\vartheta\in(0,\pi)$ be the contact angle, chosen so that $\langle N,X\rangle=\cos\vartheta$ on the unit sphere. We use $\theta$ for the angular coordinate as before. Since the surface is outside the unit ball, its outward unit conormal along $\rho=1$ is $X_\rho/\Lambda$ and
\begin{equation}\label{eq15}
X=\cos\vartheta\,N-\sin\vartheta\,\frac{X_\rho}{\Lambda}.
\end{equation}
The negative sign follows from the fact that moving into the surface, in the direction of decreasing $\rho$, increases $|X|$.

Put $T=X_\theta/\Lambda$ on the boundary, and let $\nabla$ denote the Euclidean connection on $\mathbb R^3$. By direct calculation of the tangential and normal parts, we have
\begin{align*}
\nabla_TT=-\frac{1+\Lambda_\rho/\Lambda}{\Lambda}\frac{X_\rho}{\Lambda} -\frac{A}{\Lambda^2}N.
\end{align*}
Taking its inner product with $X$ and using $\langle \nabla_TT,X\rangle=-1$, we get
\begin{equation}\label{eq16}
1+\frac{\Lambda_\rho}{\Lambda} =-\frac{\Lambda}{\sin\vartheta} +\frac{A\cos\vartheta}{\Lambda\sin\vartheta} \qquad(\rho=1).
\end{equation}
Define
\begin{align*}
\sigma=\frac{\rho^2\Lambda}{|A|} =\frac{1+|g|^2}{2|g'|},\qquad C=\frac{|A|}{\sin\vartheta}>0, \qquad D=-\operatorname{sgn}(A)\cot\vartheta.
\end{align*}
The function $\sigma$ extends positively across the puncture and across poles of $g$. Equation~\eqref{eq16} becomes
\begin{equation}\label{eq17}
\sigma_\rho=\sigma-C\sigma^2-D\qquad(\rho=1).
\end{equation}
The pullback of the spherical metric by the Gauss map is $\sigma^{-2}|dz|^2$, so
\begin{align*}
\Delta\log\sigma=\frac{1}{\sigma^2}.
\end{align*}
The assumption on the boundary regularity gives $\sigma\in C^{1,\gamma}(\overline{\mathrm D})$. Since $\sigma>0$, this is an elliptic equation with real analytic coefficients. The unit circle is real analytic, and~\eqref{eq17} expresses the outward normal derivative of $\sigma$ as a polynomial in $\sigma$. Since the boundary condition~\eqref{eq17} is also real analytic, $\sigma$ extends real analytically across $\partial\mathrm D$ by \cite[Theorem~2]{Mueller}. In particular, $\sigma$ is smooth up to the boundary.

Now, we relate the capillary boundary condition to the behavior at the end through the meromorphic structure of the Gauss map $g$. For this purpose, we consider its Schwarzian derivative. We remark that Jim\'enez  \cite{Jimenez} also used the Schwarzian derivative of the developing map to classify constant curvature metrics on annuli with constant geodesic curvature on each boundary component. Put
\begin{align*}
S(g)=\frac{g^{\prime\prime\prime}}{g'}-\frac32\left(\frac{g''}{g'}\right)^2.
\end{align*}
This is a holomorphic function on the disk, including at poles of $g$, since $S(g)=S(1/g)$. Direct differentiation gives
\begin{equation}\label{eq18}
\frac12S(g)=-\frac{\sigma_{zz}}{\sigma}.
\end{equation}
Also, on $\rho=1$,
\begin{align*}
\operatorname{Im}(z^2\sigma_{zz}) =\frac12(\sigma_\theta-\sigma_{\rho\theta}).
\end{align*}
Differentiating~\eqref{eq17} in $\theta$ and using~\eqref{eq18}, we obtain
\begin{equation}\label{eq19}
\operatorname{Im}(z^2S(g)) =\frac{\sigma_{\rho\theta}-\sigma_\theta}{\sigma} =-2C\sigma_\theta.
\end{equation}

Let $H$ be the holomorphic function whose real part is the harmonic extension of the boundary values of $\sigma$, normalized by $\operatorname{Im}H(0)=0$. On the boundary we have $\operatorname{Re}H=\sigma$ and $\sigma_\theta=-\operatorname{Im}(zH')$. Thus $z^2S(g)-2CzH'$ is holomorphic on the disk and has zero imaginary part on the boundary. Its value at the origin is zero, so
\begin{equation}\label{eq20}
z^2S(g)=2CzH'.
\end{equation}
Since $g$ has no branch point, all its poles are simple and the zeros of $1/g'$ have even order. On the simply connected disk we can choose a holomorphic function
\begin{align*}
G=(g')^{-1/2}.
\end{align*}
The singularities of $gG$ at poles of $g$ are removable. This choice of a square-root function is motivated by the work of Osserman and Schiffer \cite{OS}. We remark that $G$ and $gG$ are constant multiples of the square roots used in \cite[Lemma~3, Case~2]{OS}. In our setting, we will show that the capillary boundary condition expressed with the Schwarzian derivative implies that equality holds in the corresponding coefficient inequality.

We have
\begin{equation}\label{eq21}
\begin{gathered}
g'G^2=1,\qquad \sigma=\frac{|gG|^2+|G|^2}{2},\\
(gG)''+\frac12S(g)(gG)=0,\qquad G''+\frac12S(g)G=0.
\end{gathered}
\end{equation}
The boundary smoothness of $\sigma$ implies that $S(g)$ is smooth up to the boundary by~\eqref{eq18}. The differential equations in~\eqref{eq21} then give the same regularity for $gG$ and $G$.

Using~\eqref{eq20}, multiply the last two equations in~\eqref{eq21} by $z^2\overline{gG}$ and $z^2\overline G$, respectively, and integrate over $|z|=1$. Then
\begin{align*}
\frac1{2\pi}\int_0^{2\pi} (z^2(gG)''\overline{gG}+z^2G''\overline G)\,d\theta =-\frac{2C}{2\pi}\int_0^{2\pi}zH'\operatorname{Re}H\,d\theta.
\end{align*}
Write
\begin{align*}
gG=\sum_{m\ge0}a_mz^m,\qquad G=\sum_{m\ge0}b_mz^m,\qquad H=\sum_{n\ge0}H_nz^n.
\end{align*}
The average of $zH'H$ is zero, and the average of $zH'\overline H$ is $\sum_{n\ge1}n|H_n|^2$. Comparing the constant Fourier coefficients gives
\begin{equation}\label{eq22}
\sum_{m\ge2}m(m-1)(|a_m|^2+|b_m|^2) +C\sum_{n\ge1}n|H_n|^2=0.
\end{equation}
Since $C>0$, every term is nonnegative. It follows that $H$ is constant and $gG$ and $G$ are affine functions. Since $g(0)=0$, we have $a_0=0$ and $b_0\ne0$. The boundary value $|gG|^2+|G|^2=2\operatorname{Re}H$ is constant. Its first Fourier coefficient is $b_1\overline b_0$, which must vanish. Thus $b_1=0$, so $G$ is constant and
\begin{align*}
g(z)=az,\qquad dh=-A\frac{dz}{z},\qquad a\ne0.
\end{align*}
These are the Weierstrass data of the catenoid. The boundary circle is horizontal and centered on the axis of the catenoid. Since it lies on the unit sphere, its horizontal center is the origin, so the axis passes through the center of the sphere. This completes the proof.
\end{proof}
\begin{remark}
We note here that in case the surface has an end of multiplicity greater than 1, $f(z)$ could have a pole at the puncture. In fact, if the surface has the end asymptotic to the end of Enneper's surface, $f(z)$ has a pole of order 2.
\end{remark}

\begin{remark}
Note that Enneper's surface $X:\mathbb{C} \rightarrow \Sigma$ has the Weierstrass data $(g,dh)=(z,z dz)$ defined on $\mathbb {C}$. As \eqref{eq3} shows, the second fundamental form of the surface is given by $-\text{Re}\left(dz(v)\cdot dz(v)\right)$. Thus the images of the two perpendicular rays $z=r$ and $z=re^{i\pi/2}$ are lines of curvature. In this coordinate, we have
\begin{align*}
\Phi=-1,\qquad f(z)=-z^2,\qquad \beta=-r^2\sin 2\theta.
\end{align*}
Hence $\beta$ is zero on those rays but does not vanish identically on the surface. Here the end corresponds to $z=\infty$. In the coordinate $w=1/z$ at the end, the Hopf differential is $-w^{-4}dw^2$, which gives the pole of order 2 described in the preceding remark.
\end{remark}

\begin{remark}
In the free boundary case, similar classification results were obtained by Park and Pyo \cite{PP} and Mazet and Mendes \cite{MM}. For general capillary angles, to the author's knowledge, no such classification was previously known in the immersed setting of Theorem~\ref{thm1}. Here, by making effective use of the Weierstrass data, we obtain this classification in $\mathbb R^3$.
\end{remark}

\begin{remark}\label{lee-yeon}
The relation between the Gauss map and a Liouville equation was used in \cite{LY} to study free boundary minimal annuli. As explained in \cite[Remark~4.2]{LY}, the boundary condition depends on the parametrization as well as the shape of the Gauss image, which makes it difficult to determine the Gauss map. In our catenoidal case, the embedded end allows the Gauss map to extend without branching across the puncture, and $\sigma$ extends as a positive function. 
\end{remark}

\section{Proof of Theorem~\ref{thm2}}\label{sec4}
In this section we prove the nonexistence of such a surface when its boundary is embedded and weakly convex in $\mathrm S^2$. Note that, in the free boundary case of Theorem~\ref{thm1}, Lemma~\ref{lem1} and~\eqref{eq25} show that the geodesic curvature of the boundary in $\mathrm S^2$ is identically zero for a planar end and has a fixed nonzero sign for a catenoidal end. Thus, if the boundary is embedded, it is weakly convex, and strictly convex unless the surface is part of a plane; see \cite[Lemma~2.2]{Ghomi}. For an Enneper end, however, this argument does not give the same sign control.

\begin{proof}[Proof of Theorem~\ref{thm2}]
We begin with Hopf's arguments \cite{choe,H1}. With the notation of the previous section, we write
\begin{equation}\label{eq23}
\Phi(z)=\frac{\mathcal L-\mathcal N}{2}-i\mathcal M, \qquad f(z)=z^2\Phi(z)=\alpha+i\beta.
\end{equation}
We first calculate the curvature of the boundary in terms of $X$ and $\Lambda$. Put $\gamma(\theta)=X(e^{i\theta})$ and let $s$ be its arclength. On $\rho=1$, the free boundary condition and the fact that the surface is outside the unit ball give
\begin{align*}
X_\rho=-\Lambda X,\qquad |X_\theta|=\Lambda, \qquad N=\frac{X_\rho\times X_\theta}{\Lambda^2}.
\end{align*}
Thus
\begin{equation}\label{eq24}
\frac{d^2\gamma}{ds^2} =\frac{X_{\theta\theta}}{\Lambda^2} -\frac{\Lambda_\theta}{\Lambda^3}X_\theta =-\gamma+\kappa N,
\end{equation}
where
\begin{equation}\label{eq25}
\kappa=\frac{\langle X_{\theta\theta},N\rangle}{\Lambda^2} =\frac{\langle X_{\theta\theta},X_\rho\times X_\theta\rangle}{\Lambda^4} =-\frac{\alpha}{\Lambda^2}.
\end{equation}
Here $\kappa$ is the normal curvature of the boundary on the surface and also its geodesic curvature on $\mathrm S^2$ with respect to $N$. The last equality in~\eqref{eq25} follows from $b(\partial_\theta,\partial_\theta)=-\alpha$. Thus weak convexity of the boundary is equivalent to $\kappa$ not changing sign. This should be distinguished from the geodesic curvature of the boundary on $\Sigma$, which, with respect to the outward unit conormal, is
\begin{align*}
-\frac{\Lambda+\Lambda_\rho}{\Lambda^2}=1, \qquad\text{or equivalently}\qquad \Lambda_\rho=-\Lambda-\Lambda^2\quad(\rho=1).
\end{align*}
This identity holds for every exterior free boundary surface in these coordinates. Since the boundary is a line of curvature, $\beta=0$ there.

We now conclude with the embeddedness of the boundary. Since the geodesic curvature $\kappa$ on $\mathrm S^2$ does not change sign, the boundary curve is spherically convex and stays in a closed hemisphere; see \cite[Lemma~2.2]{Ghomi}. Hence there is a unit vector $a$ such that $\langle\gamma,a\rangle\ge0$ along the boundary. Recall that the flux vector is given by
\begin{equation}\label{eq26}
\operatorname{Flux}([\gamma])=\int_\gamma\nu(s)\,ds,
\end{equation}
where $\gamma$ is a closed oriented curve and $\nu$ is its compatible unit conormal. The flux vector depends only on the homology class of $\gamma$. 

The asymptotic condition implies that each coordinate function of $X-X_E$ extends harmonically across the puncture. Hence the flux vector around the end agrees with that of Enneper's surface and vanishes. Since $\mathrm D^\prime$ is an annulus, the boundary flux also vanishes. Since the outward unit conormal along the boundary is $\nu=X_\rho/\Lambda=-\gamma$, equation~\eqref{eq26} gives
\begin{align*}
0=-\left\langle\operatorname{Flux}([\partial\Sigma]),a\right\rangle =\int_{\partial\Sigma}\langle\gamma,a\rangle\,ds.
\end{align*}
The integrand is nonnegative, so the boundary curve lies in the bounding great circle. It follows that $\kappa=0$ and therefore $\alpha=\beta=0$ on the boundary. The function $f=z^2\Phi$ is holomorphic near the boundary from inside the disk and vanishes on $|z|=1$. Schwarz reflection extends $f$ holomorphically across the circle. The identity theorem then gives $f\equiv0$ on $\mathrm D^\prime$, and hence $\Phi\equiv0$. The surface is a part of the plane, which contradicts the assumption on the end. This completes the proof.
\end{proof}

\begin{remark}
The weak convexity assumption excludes boundary umbilic points of rotation index $-1/4$. Indeed, by~\eqref{eq25}, $\alpha$ does not change sign along the boundary. The function $f=z^2\Phi$, which is real on $|z|=1$, extends holomorphically across the circle by Schwarz reflection. At an isolated boundary umbilic point its order of vanishing must therefore be even, say $2m$. Its boundary rotation index is $-m/2$; see \cite[Lemma~2]{choe}. If $f$ vanishes on a boundary arc, analytic continuation gives $\Phi\equiv0$, and the surface is planar.
\end{remark}

\begin{remark}\label{rmk:enneper}
We expect that there exist annular exterior free boundary minimal immersions with embedded boundary, finite total curvature and one Enneper end, when the weak convexity assumption in Theorem~\ref{thm2} is removed. Any such example must have exactly four boundary umbilic points of rotation index $-1/4$ and no interior umbilic points. Indeed, the Enneper end gives $f=z^2\Phi$ a pole of order 2 at the puncture. Since $f$ is real on $|z|=1$, Schwarz reflection shows that $z^2f$ is a polynomial of degree 4. On the other hand, zero flux gives $\int_{\partial\Sigma}\gamma\,ds=0$, so the convex hull of the embedded spherical boundary contains the origin. By \cite[Theorem~1.2]{Ghomi}, the boundary has at least four points where its geodesic curvature on $\mathrm S^2$ vanishes. These are zeros of $f$ by~\eqref{eq25}, so there are exactly four, all simple, and no interior zeros. The existence of such a surface remains open, and we leave this problem for further study.
\end{remark}

\begin{remark}
The same argument applies to higher order Enneper ends, with the same boundary regularity, free boundary, embedded boundary, zero flux and weak convexity assumptions. In fact, once $f$ vanishes on the boundary, Schwarz reflection and analytic continuation give $\Phi\equiv0$, without using its behavior at the puncture. Thus, in the same annular setting, the planarity conclusion does not require finite total curvature or a prescribed asymptotic behavior at the end.
\end{remark}

\begin{remark}\label{rmk:enneper-asymptotics}
If the asymptotic condition~\eqref{eq5} is weakened to the assumption that the end has multiplicity 3, its flux vector need not vanish. The ends of the minimal annulus given by~\eqref{eq6} satisfy this weaker condition. Indeed, writing $a=(R^2+R^{-2})^{-1}$, its height differential near $z=0$ is $dh=(-az^{-3}+z^{-1}-az)\,dz$. Thus the end has nonzero flux \cite[Section~2.3.10]{K}.
\end{remark}

\section{Proof of Theorem~\ref{thm3}}\label{sec5}
As in \cite{choe}, we calculate Hopf's rotation index at the boundary umbilic points and vertices of the surface to prove Theorem~\ref{thm3}.

\begin{proof}[Proof of Theorem~\ref{thm3}]

If the surface is totally umbilic, it must be part of the plane. Assume otherwise. Suppose that a vertex $p$ has interior angle $\xi<\pi$, and let $P=T_p\Sigma$. Since both boundary arcs lie in $\overline U$, the exterior condition implies that at least one defining face through $p$ is tangent to $P$. Otherwise, the surface would enter $U$ along an interior direction of the smaller tangent sector. Such a face cannot be spherical, since the normal curvatures in the two independent principal directions would have the same nonzero sign, contradicting minimality.

Thus all defining faces tangent to $P$ are planar. Choose a unit normal $n$ to $P$ such that $U\subset\{h\le0\}$, where $h(x)=\langle x-p,n\rangle$. If $h<0$ at a point of either adjacent boundary arc sufficiently close to $p$, since $\xi<\pi$, moving slightly into the interior of the surface from such a point would give a point inside every defining ball and half-space. This contradicts the exterior condition, so both arcs lie in $P$. Choose the orientation so that $N(p)=n$. Since each arc is a line of curvature contained in $P$, $\langle N,n\rangle$ is constant along each arc, and hence $N=n$. Minimality then implies that the Hopf differential vanishes there, and Schwarz reflection and the identity theorem imply that the surface is planar, a contradiction. Thus $\xi\ge\pi$. Since the two boundary tangent lines are distinct, $\xi=\pi$ is excluded, and hence $\xi>\pi$.

The differential equation of the lines of curvature is given by
\begin{align*}
\mathcal M\,du^2+(\mathcal N-\mathcal L)\,du\,dv-\mathcal M\,dv^2=0,
\end{align*}
which is equivalent to $\operatorname{Im}(\Phi\,dz^2)=0$. At an isolated interior umbilic point $p$, the rotation index is given by \cite{H1}
\begin{equation}\label{eq27}
I_{\mathrm{rot}}(p)=-\frac{1}{4\pi}\delta(\arg\Phi)=-\frac m2
\end{equation}
when $\Phi$ has a zero of order $m$. Here $\delta(\arg\Phi)$ denotes the total change of $\arg\Phi$ along a small positively oriented loop enclosing $p$ and no other zero. A pole of order $m$ has index $m/2$.

Let us look at the order of $\Phi(z)$ at the puncture. As in the proof of Lemma~\ref{lem1}, the Weierstrass data in~\eqref{eq11} and~\eqref{eq12} give
\begin{align}\label{eq28}
\Phi(z)=-\frac{dh}{g\,dz}g'(z)=-n c_n d_{-2}z^{n-3}+O(z^{n-2}),
\end{align}
where $c_n d_{-2}\ne0$.

Assume that the surface has a catenoidal end. Since $n=1$, equation~\eqref{eq28} shows that $\Phi$ has a pole of order 2 at the puncture. Thus
\begin{align}\label{eq29}
I_{\text{rot}}(z=0) =1 .
\end{align}
As in the proof of Hopf's theorem \cite{H1}, at an isolated interior umbilic point $q$,
\begin{align}\label{eq30}
I_{\text{rot}}(q) \le-\frac{1}{2}
\end{align}
since  $\Phi(z)$ has a zero at $z=q$.

Each regular boundary arc is a line of curvature by the Terquem-Joachimsthal theorem, since it meets a sphere or a plane at a constant contact angle. Let $t \in \partial \Sigma$ denote a boundary umbilic point or a vertex point of the boundary. By the angle estimate above and the estimation of the rotation index of the boundary umbilic points \cite[Lemma~2]{choe}, we have
\begin{align}\label{eq31}
I_{\text{rot}}(t) \le-\frac{1}{4}.
\end{align}
At a regular boundary point which is not umbilic, the boundary index is zero.

Since $\Phi$ is not identically zero, its interior zeros are isolated. Reflection across the regular boundary arcs and the local meromorphic extension at the vertices \cite[Lemmas~2 and~3]{choe}, together with the expansion at the puncture, exclude accumulation of singularities. Hence the curvature-line field has only finitely many singularities on the compactified parameter disk.

After filling in the puncture, the parameter domain is $\overline{\mathrm D}=\{z\in\mathbb C\mid |z|\le1\}$. Applying \eqref{eq7} to the reflected line field on the double of this disk and dividing by two, we have
\begin{align*}
\sum_p I_{\mathrm{rot}}(p)=\chi(\overline{\mathrm D})=1,
\end{align*}
where the puncture is included in the sum. Equations \eqref{eq29} and \eqref{eq30}, together with \eqref{eq31}, imply that if $\partial\Sigma$ has at least one vertex,
\begin{align*}
\sum_p I_{\mathrm{rot}}(p)\le1-\frac14<1.
\end{align*}
This is a contradiction. Thus it follows that the surface with a catenoidal end cannot exist in this case.

Now assume that the surface has a planar end. Since $n>1$, \eqref{eq28} implies that $\big{|}z\Phi(z)\big{|} \le C$ for some fixed constant $C$. We see that
\begin{align*}
I_{\text{rot}}(z=0) \le \frac{1}{2} .
\end{align*}
In the same way, all the other indices are nonpositive, so their sum with the index at the puncture is at most $1/2<1$. Thus the surface cannot exist unless it is totally umbilic, and hence it must be part of the plane. This completes the proof.
\end{proof}

\begin{remark}
The case in which the boundary is smooth and passes from one face to another, without being contained in a single face, is not covered by the nonexistence assertion for catenoidal ends in Theorem~\ref{thm3}. At a smooth transition point the boundary angle is $\pi$, so the strictly negative vertex index used in~\eqref{eq31} is no longer guaranteed. Analytic continuation across such a point would not by itself provide the negative index, which is crucial for the above argument. The regularity of the immersion across a change of face is a separate question, which we leave for further study.
\end{remark}

\end{document}